

  \newcount\fontset
  \fontset=1
  \def \dualfont#1#2#3{\font#1=\ifnum\fontset=1 #2\else#3\fi}

  \dualfont\bbfive{bbm5}{cmbx5}
  \dualfont\bbseven{bbm7}{cmbx7}
  \dualfont\bbten{bbm10}{cmbx10}

  \font \eightbf = cmbx8
  \font \eighti = cmmi8 \skewchar \eighti = '177
  \font \eightit = cmti8
  \font \eightrm = cmr8
  \font \eightsl = cmsl8
  \font \eightsy = cmsy8 \skewchar \eightsy = '60
  \font \eighttt = cmtt8 \hyphenchar\eighttt = -1

  \font \sixi = cmmi6 \skewchar \sixi = '177
  \font \sixrm = cmr6
  \font \sixsy = cmsy6 \skewchar \sixsy = '60
  \font \tensc = cmcsc10
  
  \font \titlefont = cmbx12
  \scriptfont \bffam = \bbseven
  \scriptscriptfont \bffam = \bbfive
  \textfont \bffam = \bbten

  \newskip \ttglue

  \def \eightpoint {\def \rm {\fam0 \eightrm }%
  \textfont0 = \eightrm
  \scriptfont0 = \sixrm \scriptscriptfont0 = \fiverm
  \textfont1 = \eighti
  \scriptfont1 = \sixi \scriptscriptfont1 = \fivei
  \textfont2 = \eightsy
  \scriptfont2 = \sixsy \scriptscriptfont2 = \fivesy
  \textfont3 = \tenex
  \scriptfont3 = \tenex \scriptscriptfont3 = \tenex
  \def \it {\fam \itfam \eightit }%
  \textfont \itfam = \eightit
  \def \sl {\fam \slfam \eightsl }%
  \textfont \slfam = \eightsl
  \def \bf {\fam \bffam \eightbf }%
  \textfont \bffam = \bbseven
  \scriptfont \bffam = \bbfive
  \scriptscriptfont \bffam = \bbfive
  \def \tt {\fam \ttfam \eighttt }%
  \textfont \ttfam = \eighttt
  \tt \ttglue = .5em plus.25em minus.15em
  \normalbaselineskip = 9pt
  \def \MF {{\manual opqr}\-{\manual stuq}}%
  \let \sc = \sixrm
  \let \big = \eightbig
  \setbox \strutbox = \hbox {\vrule height7pt depth2pt width0pt}%
  \normalbaselines \rm }



  \newcount \secno \secno = 0
  \newcount \stno \stno = 0
  \newcount \eqcntr \eqcntr= 0

  \def \ifn #1{\expandafter \ifx \csname #1\endcsname \relax }

  \def \track #1#2#3{\ifn{#1}\else {\tt\ [#2 \string #3] }\fi}

  \def \advseqnumbering {\global \advance \stno by 1 \global \eqcntr=0}

  \def \current {\number \secno \ifnum \number \stno = 0 \else
    .\number \stno \fi }

  \def \laberr#1#2{\message{*** RELABEL CHECKED FALSE for #1 ***}
      RELABEL CHECKED FALSE FOR #1, EXITING.
      \end}

  \def \syslabel#1#2{%
    \ifn {#1}%
      \global \expandafter 
      \edef \csname #1\endcsname {#2}%
    \else
      \edef\aux{\expandafter\csname #1\endcsname}%
      \edef\bux{#2}%
      \ifx \aux \bux \else \laberr{#1=(\aux)=(\bux)} \fi
      \fi
    \track{showlabel}{*}{#1}}

  \def \subeqmark #1 {\global \advance\eqcntr by 1
    \edef\aux{\current.\number\eqcntr}
    \eqno {(\aux)}
    \syslabel{#1}{\aux}}

  \def \eqmark #1 {\advseqnumbering
    \eqno {(\current)}\syslabel{#1}{\current}}

  \def \label #1 {\syslabel{#1}{\current}}

  \def \lcite #1{(#1\track{showcit}{$\bullet$}{#1})}

  \def \cite #1{[{\bf #1}\track{showref}{\#}{#1}]}

  \def \scite #1#2{{\rm [\bf #1\track{showref}{\#}{#1}{\rm \hskip 0.7pt:\hskip 2pt #2}\rm]}}


 \def \Headlines #1#2{\nopagenumbers
    \advance \voffset by 2\baselineskip
    \advance \vsize by -\voffset
    \headline {\ifnum \pageno = 1 \hfil
    \else \ifodd \pageno \tensc \hfil \lcase {#1} \hfil \folio
    \else \tensc \folio \hfil \lcase {#2} \hfil
    \fi \fi }}

  \def \Title{\centerline{\titlefont \ucase{\titletextOne}}
    \edef\aux{\titletextTwo}\edef\bux{\null}%
    \ifx\aux\bux \else 
      \smallskip
      \centerline{\titlefont \ucase{\titletextTwo}}
      \fi}

  \def \Date #1 {\footnote {}{\eightit Date: #1.}}

  \long \def \Abstract #1{\begingroup
  \bigskip \bigskip \noindent
  {ABSTRACT.} #1\par \par \endgroup }


  \def \ucase #1{\edef \auxvar {\uppercase {#1}}\auxvar }
  \def \lcase #1{\edef \auxvar {\lowercase {#1}}\auxvar }

  \def \section #1{\global\def \SectionName{#1}\stno = 0 \global
\advance \secno by 1 \bigskip \bigskip \goodbreak \noindent {\bf
\number \secno .\enspace #1.}\medskip \noindent \ignorespaces}

  \long \def \sysstate #1#2#3{%
    \advseqnumbering
    \medbreak \noindent 
    {\bf \current.\enspace #1.\enspace }{#2#3\vskip 0pt}\medbreak }
  \def \state #1 #2\par {\sysstate {#1}{\sl }{#2}}
  \def \definition #1\par {\sysstate {Definition}{\rm }{#1}}
  \def \remark #1\par {\sysstate {Remark}{\rm }{#1}}


  \def \proof {\medbreak \noindent {\it Proof.\enspace }}
  \def \proofend {\ifmmode \eqno \quadrado \else \hfill \quadrado
\looseness = -1 \medbreak \fi }

  \def \$#1{#1 $$$$ #1}
  \def \=#1{\buildrel \hbox{\sixrm #1} \over =}

  \def \Item #1{\smallskip \item {{\rm #1}}}
  \newcount \zitemno \zitemno = 0

  \def \zitemplus {\global \advance \zitemno by 1 \relax}
  \def \rzitem{\romannumeral \zitemno}

  \def \zitemmark #1 {\syslabel{#1}{\rzitem}}

  \newcount \nitemno \nitemno = 0
  
  \def \nitem {\global \advance \nitemno by 1 \Item {{\rm(\number\nitemno)}}}

  \newcount \aitemno \aitemno = -1
  \def \boxlet#1{\hbox to 6.5pt{\hfill #1\hfill}}
  
  \def \aitemconv{\ifcase \aitemno a\or b\or c\or d\or e\or f\or g\or
h\or i\or j\or k\or l\or m\or n\or o\or p\or q\or r\or s\or t\or u\or
v\or w\or x\or y\or z\else zzz\fi}
  \def \aitem {\global \advance \aitemno by 1\Item {(\boxlet \aitemconv)}}
  \def \aitemmark #1 {\syslabel{#1}{\aitemconv}}

  \newcount \footno \footno = 1
  \newcount \halffootno \footno = 1
  \def \footcntr {\global \advance \footno by 1
  \halffootno =\footno
  \divide \halffootno by 2
  $^{\number\halffootno}$}
  \def \fn#1{\footnote{\footcntr}{\eightpoint#1\par}}

  \begingroup
  \catcode `\@=11
  \global\def\eqmatrix#1{\null\,\vcenter{\normalbaselines\m@th%
      \ialign{\hfil$##$\hfil&&\kern 5pt \hfil$##$\hfil\crcr%
	\mathstrut\crcr\noalign{\kern-\baselineskip}%
	#1\crcr\mathstrut\crcr\noalign{\kern-\baselineskip}}}\,}
  \endgroup


  \font\mf=cmex10
  \def\union {\mathop{\raise 9pt \hbox{\mf S}}\limits}
  \def\inters{\mathop{\raise 9pt \hbox{\mf T}}\limits}

  \def \C {{\bf C}}
  \def \<{\left \langle \vrule width 0pt depth 0pt height 8pt }
  \def \>{\right \rangle }

  \def \and {\hbox {,\quad and \quad }}
  \def \calcat #1{\,{\vrule height8pt depth4pt}_{\,#1}}

  \def \for #1{,\quad \forall\,#1}
  \def \quadrado {\hbox {$\sqcap \!\!\!\!\sqcup $}}
  
  \def \stress #1{{\it #1}\/}
  \def \inv {^{-1}}
  \def \*{\otimes}

  \newcount \bibno \bibno = 0
  \def \newbib #1{\global\advance\bibno by 1 \edef #1{\number\bibno}}

  \def \bibitem #1#2#3#4{\smallskip \item {[#1]} #2, ``#3'', #4.}

  \def \references {
    \begingroup
    \bigskip \bigskip \goodbreak
    \eightpoint
    \centerline {\tensc References}
    \nobreak \medskip \frenchspacing }



\def \caldef #1{\global \expandafter \edef \csname #1\endcsname {{\cal #1}}}

  \input amssym  

  \font \fivesl = cmsl5  

  \newbib{\Actions} 
  \newbib{\infinoa} 
  \newbib{\Skand} 
  \newbib{\Lawson} 
  \newbib{\Paterson} 
  \newbib{\Renault} 

  \def\c{\subseteq}
  \def\F{{\cal F}}
  \def\G{{\cal G}}
  \def\H{{\cal H}}

  \def \Bitem {\Item {$\bullet $}}
  \def\inters{\raise 0.4pt \hbox{$\cap\kern0.5pt$}}
  \def\tightbox{{\hbox{\fivesl tight}}}                     
  
  \def\b{\beta}
  \def\g{\gamma}
  \def\Gz{\G^{(0)}}
  \def\ci#1{\scite{\Actions}{#1}}
  \def\d{d}   \def\r{r}
  \def\its{\Cap}
  \def\piu{\pi_u}
  \def\t{\theta}
  
  \def\actspec{\theta}

  \def\tb{{\hbox{\fivesl tight}}}                     
  \font\bigcal=rsfs10
  \font\smallcal=rsfs8 
  \def\S{\hbox{\bigcal S}\kern 2.5pt}
  \def\s{\hbox{\smallcal S}\,}
  \def\E{\hbox{\bigcal E}\,}
  \def\e{\hbox{\smallcal E}\,}
  \def\Eh{\, \widehat {\!\E}}
  \def\Et{\Eh_\tb}
  \def\D{{\cal D}}
  \def\Cstght{C^*_\tb(\S)}

  \def\Bas{{\cal B}}
  \def\Bast{\widehat \Bas_\tb}
  \def\X{{\cal X}}
  \def\slice{slice}

\long\def\Abstract#1{\midinsert \narrower \narrower \eightpoint
\noindent #1 \endinsert}


  \Headlines {TOTALLY DISCONNECTED GROUPOIDS AND AMPLE SEMIGROUPS} {R.~Exel}
  \font \titlefont = cmbx9
  \null\vskip -1cm
  \centerline{\titlefont RECONSTRUCTING A TOTALLY DISCONNECTED}
  \smallskip
  \centerline{\titlefont  GROUPOID FROM ITS AMPLE SEMIGROUP}

  \footnote{\null} 
  {\eightrm 2000 \eightsl Mathematics Subject Classification:
  \eightrm 
  22A22, 
  20M18, 
  20M30, 
  46L55. 
  }

  \bigskip
  \centerline{\tensc 
    R.~Exel\footnote{*}{\eightpoint Partially supported by
CNPq.}}

  \bigskip
  \Date{25 Jul 2009}

\Abstract{We show that a (not necessarily Hausdorff) \'etale, second
countable groupoid $\G$ with totally disconnected unit space may be
reconstructed solely from the algebraic structure of its ample
semigroup $\s$.  We also show that $C^*(\G)$ possesses a universal property
related to tight representations of $\s$.}

\section {Introduction}
  Recall that an open subset $S$ of a (not necessarily Hausdorff)
\'etale groupoid $\G$ is said to be a \stress{\slice} (sometimes also a
\stress{$\G$-set}, or a \stress{bissection}) if the domain and range maps are injective on
$S$.  The collection of all {\slice}s forms an inverse semigroup which
has often been studied alongside the groupoid itself.
  Among other things the natural action of this semigroup on the unit space
$\Gz$ highlights the dynamical nature of groupoids.

The set $\S$ formed by all compact (open) slices is 
called the \stress{ample semigroup} of $\G$ \scite{\Renault}{2.10}.
A compact {\slice} $S$ is somewhat special because its
characteristic function $1_S$ is in $C_c(G)$ and hence also in
$C^*(\G)$.  
We
therefore get a map
  $$
  \rho:\S \to C^*(\G),
  \eqmark TheMapIntroduced
  $$
  defined by $\rho(S) = 1_S$, which is multiplicative and satisfies
$\rho(S\inv) = \rho(S)^*$.  In other words, $\rho$ is a \stress{representation} of
$\S$ in $C^*(\G)$.

From now on we will specialize to the situation in which $\Gz$ is
totally disconnected.  In this case one may prove that $\G$ contains
many compact {\slice}s, enough to provide a basis for the topology of
$\G$.  The starting point for the present investigation was our
interest in determining precisely how much information does $\S$
contain about $\G$ and, in particular, one may ask the
following:

\sysstate{Question}{\rm}{\label Questao Is it possible to reconstruct
$\G$ solely from the algebraic structure of the inverse semigroup $\S$ of compact {\slice}s?}

The literature contains several constructions of groupoids associated
to inverse semigroups, such as the \stress{associated groupoid}
described by Lawson in \scite{\Lawson}{Section 3.1}, Kellendonk's
topological groupoid \scite{\Lawson}{Section 9.2}, or Paterson's
groupoid of germs for the action described in \scite{\Paterson}{Proposition 4.3.2}
  but none of these give the expected result when the ample
semigroup is fed as input.

  The purpose of this note is to give an affirmative answer to the
above question, assuming $\G$ to be second countable, by using the
theory of \stress{tight representations} of semilattices and
inverse semigroups introduced by the author in \cite{\Actions}.  
  En passant we prove that the unitization of $C^*(\G)$ is the
universal C*-algebra for tight representations of $\S$ in the sense
that these factor through the representation $\rho$ of
\lcite{\TheMapIntroduced}.

Question \lcite{\Questao} has already been addressed by Renault (see
the paragraph immediately following \scite{\Renault}{2.10}) but under
the assumption that the unit space of $\G$ is known.  Without this
crucial piece of information one is faced with the task of
reconstructing the unit space of $\G$ from $\S$.  Our strategy is to
recover $\Gz$ from the idempotent semilattice $\E$ of $\S$ by means of
its \stress{tight spectrum}.  The action of $\S$ on $\Gz$ is then
easily recovered and $\G$ pops out as the groupoid of germs for this
action, as long as we use the appropriate notion of germs
\scite{\Paterson}{page 140}.

I would like to acknowledge an interesting conversation with Alex
Kumjian in which the above question was raised,  and for which thanks
are due!

\section {Preliminaries}
  \label PrelimSection
  In this section we shall briefly survey the main results from
\cite{\Actions} that we need here.  We will begin with some basic
facts about inverse semigroups, all of which may be
found in greater detail in many books on the subject such as
\cite{\Paterson} or \cite{\Lawson}.
   References for groupoids may also be found in \cite{\Renault} in
the Hausdorff case, in \cite{\Skand} in the general non-Hausdorff
locally compact case, and in the first sections of \cite{\Actions}
which is dedicated to the non-Hausdorff \'etale case,  
the situation which will concern us here.

  Recall that a semigroup $\S$ 
  is said to be an \stress{inverse semigroup} if for every $s\in \S$
there exists a unique $s^*\in \S$ such that
  $ss^*s=s$ and   $s^*ss^*=s^*$.
An element $0\in \S$ such that $0s = s\kern0.3pt 0 = 0$, for all $s\in \S$, is said
to be a \stress{zero element}.  All of our inverse semigroups will
henceforth be assumed to contain a zero element.

An element $e\in \S$ is said to be \stress{idempotent} if $e^2=e$.  In
this case $e$ is necessarily self-adjoint, meaning that $e^*=e$.  The
subset $\E\c \S$ consisting of all idempotent elements is a commutative
subsemigroup.  It is also a semilattice under the order relation
  $$
  e\leq f \iff ef = e
  \for e,f\in \E.
  $$

Given $e,f\in \E$, we say that $e$ is \stress{orthogonal} to $f$, in
symbols $e\perp f$, when $ef = 0$.  Otherwise we say that $e$
\stress{intersects} $f$, writing $e\its f$.
Given subsets $Z, F\c \E$, we will say that $Z$ is a \stress{cover} \ci{11.5}
for
$F$, if $Z\c F$ and for every $f\in F$ there exists $z\in Z$
such that $z\its f$.

A \stress{representation} of $\S$ in a unital C*-algebra $A$ is a map
  $$
  \pi:\S \to A,
  $$
  such that $\pi(st)=\pi(s)\pi(t)$, and $\pi(s^*)=\pi(s)^*$, for all
$s,t\in \S$, and such that $\pi(0)=0$.  We shall moreover say that
$\pi$ is \stress{tight} \ci{11.6,  13.1}
if, whenever one is given finite subsets
$X,Y\c \E$, and a finite cover $Z$ for
  $$
  \E^{X,Y} := \{f\in \E: f\leq x,\ \forall x\in X, \hbox{ and } f\perp y,\ \forall y\in Y\}, 
  $$
  one has that
  $$
  \bigvee_{z\in Z}   \pi(z) =
  \prod_{x\in X} \pi(x)\prod_{y\in Y} \big(1-\pi(y)\big).
  \eqmark TightCondition
  $$
  Here ``$\vee$" refers to the operation of taking the supremum of a
commuting set of projections (notice that the commutativity of $\E$
implies that the $\pi(z)$ commute with each other) which, in case of
two projections $p$ and $q$, is defined by
  $$
  p\vee q = p+q-pq.
  $$

The following is a central concept,  implicit in \cite{\Actions}:

\definition \label DefCStarTight The \stress{tight} C*-algebra of $\S$,
denoted $\Cstght$, is the universal unital C*-algebra with one
generator for each element of $\S$, subject to the relations
stating that the standard map
  $$
  \piu: \S \to \Cstght
  $$
  is a tight representation.  Henceforth $\piu$ will be referred to as
the \stress{universal tight representation} of $\S$.

The requirement that $\Cstght$ be \stress{unital} is necessary because
the right hand side of \lcite{\TightCondition} explicitly mentions
``1".  However, in case $\E$ admits no finite cover,
\lcite{\TightCondition} will never occur with $X\neq\emptyset$, and
hence the right hand side of \lcite{\TightCondition}, once expanded
out, will actually \stress{not} involve ``1".  This phenomena has
already been observed in \scite{\infinoa}{Section 8}.

As is the case with universal objects, 
  $\pi_u$ possesses the following universal property:  given any tight
representation $\pi$ of $\S$ in a unital C*-algebra $A$,  there exists
a unique unital *-homomorphism $\phi: \Cstght \to A$, such
that
  $\phi\circ \piu = \pi$.

In the remainder of this section we will briefly recall the
description of $\Cstght$ as a groupoid C*-algebra given in
\cite{\Actions}.

A representation of $\E$ in the algebra of complex numbers, say
  $$
  \phi:\E\to\C
  $$
  is called a \stress{character} (without the assumption that
$\phi(0)=0$, these are sometimes also called semicharacters).  Since
the elements of $\E$ are idempotent, the range of $\phi$ is contained
in the set $\{0,1\}$, so $\phi$ is necessarily the characteristic
function of some subset $\xi\c \E$.  One readily verifies that
$0\notin\xi$, and that $\xi$ is downwards directed as well as upwards
hereditary.
  Such a $\xi$ is called a \stress{filter} \ci{12.1}, and it is
easy to see that the characteristic function of any filter is a
character.  So there is a one-to-one correspondence between the set of
all characters and the set of all filters, according to which we often
identify one another.

An \stress{ultra-filter} is a filter which is not properly contained
in any filter, and we call the corresponding character an
\stress{ultra-character}.  On the other hand,  if a filter is associated to a tight
character, we call it a \stress{tight filter}.

With the topology inherited from the product topology on
$\{0,1\}^{\e}$, the set $\Eh$ of nonzero characters is a locally
compact space.  The set of all tight characters is denoted $\Et$, and
it happens to be precisely the closure \ci{12.9} of the set of all
ultra-characters within $\Eh$.  We refer to $\Et$ as the \stress{tight
spectrum} of $\E$.

If $s\in \S$ and if $\phi\in\Eh$ is such that $\phi(s^*s)\neq0$, then
the map
  $$
  \t_s(\phi): e\in \E \mapsto \phi(s^*es) \in \{0,1\}
  \eqmark DefineTheta
  $$
  is a character satisfying $\t_s(\phi)\calcat{ss^*}\neq0$ \ci{10.3}.  

Denoting by $\D_e$ the set of all characters not vanishing on a given
idempotent $e$, one may prove that
  $\t_s$ defines a homeomorphism from $\D_{s^*s}$ to $\D_{ss^*}$.
  Collectively the family of {partial homeomorphisms} $\t=
\{\t_s\}_{s\in \s}$ is an \stress{action} \ci{4.3} of $\S$ on $\Eh$,
which leaves $\Et$ invariant, in the sense that
$\t_s(\D_{s^*s}\cap\Et)\c \Et$ \ci{12.11}.  So we get, by restriction,
an action of $\S$ on $\Et$.

For our purposes the relevant notion of \stress{germs}
\cite{\Paterson}, \ci{4.17} is defined as follows: if $s_i\in \S$, and
if $\xi\in \D_{s_i^*s_i}$, for $i=1, 2$, we say that $s_1$ and $s_2$
\stress{have the same germ}\fn{Except for the case of topologically
free actions, 
this notion of germs is not the same as
the one which emphasizes an open set containing $\xi$ where $\t_{s_1}$
and $\t_{s_2}$ agree.} at $\xi$, if there exists $e\in \E$, such that
$\xi\in \D_e$, and $s_1e=s_2e$.
  If $\xi\in \D_s$, then the equivalence class for the relation
``{having the same germ}" is called the \stress{germ} of $s$ at $\xi$
and is denoted $[s,\xi]$.  The groupoid of all germs for the action of
$\S$ on $\Et$ is denoted $\G_\tb = \G_\tb(\S)$.  See \cite{\Actions} for more details.

  \def\unit#1{#1^{\null^\sim}}

The following is but a reinterpretation of \ci{13.3}.  In it we will
denote by $\unit A$ the unitization of a C*-algebra $A$, with the
convention that $\unit A = A$, in case $A$ already has a unit.

\state Theorem \label Reinterpret Given a countable inverse semigroup
$\S$, there is a *-isomorphism
  $$
  \Lambda: \Cstght \to \unit{C^*(\G_\tb)},
  $$
  such that,  for every $s\in\S$,  
  $$
  \Lambda(\piu(s)) = 1_{\X_s},
  \subeqmark LambdaPiu
  $$
  where the right hand side refers to the characteristic function of
the compact {\slice}
  $$
  \X_s = \{[s,\xi]: \xi\in \D_{s^*s}\},
  $$
  seen as an element of $C_c(\G_\tb)\c C^*(\G_\tb)$.  In addition, the
closed *-subalgebra of $\unit{C^*(\G_\tb)}$ generated by the range of
$\Lambda\circ \piu$ coincides with $C^*(\G_\tb)$.

\proof
  Initially notice that $\X_s$ is a {\slice} of $\G_\tb$ by
\ci{4.18}, which is compact by \ci{4.15}, since $\D_{s^*s}$ is
compact,  as argued in the second paragraph after \ci{10.2}.
  Therefore the characteristic function of $\X_s$ is in $C_c(\G_\tb)\c
C^*(\G_\tb)$.

We next claim that the set $\{1_{\X_s}: s\in \S\}$ generates
$C^*(\G_\tb)$ as a C*-algebra.  Denoting by $A$ the closed
*-subalgebra generated by this set, notice that, in view of \ci{3.10},
the claim will follow once we prove that $C_c(\X_s)\c A$, for all
$s\in \S$.
  
  If $f\in C_c(\X_s)$, let $g$ be the unique element of
$C_c(\X_{s^*s})$, such that $g(\d(\gamma)) = f(\gamma)$, for all
$\gamma\in\X_s$, where $\d$ refers to the source map.  Since $f$
coincides with the product of $1_{\X_s}$ and $g$ in $C^*(\G_\tb)$, our
claim will be proved once we show that $g\in A$.

Using \ci{4.16} and the Stone--Weierstrass Theorem we have that  $C_0(\Gz_\tb)$ is generated by
the set $\{1_{\X_e}: e\in \E\}$, so 
  $$
  g\in C_c(\X_{s^*s}) \c C_0(\G_\tb^{(0)}) \c A.
  $$
  This proves our claim.

  Viewing $\Cstght$ as an algebra of operators on a Hilbert space $H$
via some faithful representation, $\piu$ may be regarded as a tight
representation of $\S$ on $H$.  Applying \ci{13.3} we deduce that
there is a *-representation $\psi$ of $C^*(\G_\tb)$ on $H$ such that
  $$
  \psi(1_{\X_s}) = \piu(s)
  \for s\in \S,
  \subeqmark psiPiu
  $$
  (notice that the argument of $\rho$ in \ci{10.15} is precisely
$1_{\X_s}$ in the notation adopted here).
  By the claim we deduce that the range of $\psi$ is contained in
$C^*_\tb(\S)$, so we may view $\psi$ as a *-homomorphism
  $$
  \psi:C^*(\G_\tb) \to \Cstght.
  $$
  We next want to extend this to unitizations but, while $\Cstght$
is unital by definition, $C^*(\G_\tb)$ may or may not be unital.
However, we claim that if $C^*(\G_\tb)$ is unital then
  $$
  \psi(1) = 1.
  $$
  Assuming that $C^*(\G_\tb)$ is unital, it follows that $\Gz_\tb$ is
compact
  and hence, by \ci{4.16}, $\Et$ is also compact.  Since
$\{\D_e\}_{e\in\e}$ is a cover for $\Et$, we may extract a finite
subcover, say
  $$
  \Et = \D_{e_1} \cup \ldots \cup\D_{e_n}.
  $$
  The unit of $C^*(\G_\tb)$ is then given by
  $$
  1 = \bigvee_{k=1}^n 1_{\X_{e_k}}.
  \subeqmark UnitOfGpd
  $$

  Given any nonzero $e\in\E$, choose an ultra-filter $\xi$ containing
$e$, as done in \ci{Section 12}, and observe that $\xi\in\Et$, by
\ci{12.7}.  Therefore $\xi\in \D_{e_k}$, for some $k$, and hence
$e_k\in\xi$,  so that $e \its e_k$.
   In other words, $\{e_1,\ldots,e_n\}$ is a cover for $\E$, in the
sense of \ci{11.5}.  

Writing $\E = \E^{X, Y}$, where $X =\emptyset$,
and $Y = \{0\}$, and observing that $\piu$ is tight by definition,
we deduce that
  $$
  \bigvee_{k=1}^n \piu(e_k) = \prod_{x\in X} \piu(x)\prod_{y\in Y}
\big(1-\piu(y)\big) = 1.
  $$
  So, 
  $$
  \psi(1) \={(\UnitOfGpd)}
  \bigvee_{k=1}^n \psi(1_{\X_{e_k}}) \={(\psiPiu)}
  \bigvee_{k=1}^n \piu(e_k) =
  1.
  $$

Being a unital algebra, $C^*(\G_\tb)$ 
coincides with its unitization and $\psi$ may be viewed as a unital
*-homomorphism
  $$
  \widetilde\psi : \unit{C^*(\G_\tb)} \to \Cstght.
  $$
  On the other hand, should $C^*(\G_\tb)$ lack a unit, the natural
extension of $\psi$ to its unitization again provides a unital
homomorphism $\widetilde\psi$ as above.  Thus $\widetilde\psi$ is
available whatever the case may be.
  
Changing tack, consider the map  
  $$
  \rho: s\in \S   \mapsto 1_{\X_s} \in \unit{C^*(\G_\tb)}.
  $$
  Either using \ci{13.3} or by direct computation one proves that $\rho$
is a tight representation of $\S$, in which case the universal
property of $\Cstght$ intervenes providing a unital *-homomorphism
  $$
  \Lambda: \Cstght \to \unit{C^*(\G_\tb)},
  $$
  such that
  $$
  \Lambda(\piu(s)) = \rho(s) = 1_{\X_s}
  \for s\in\S.
  \subeqmark PiuPsi
  $$
  Noticing that $\Cstght$ is generated by $\{\piu(s): s\in \S\} \cup
\{1\}$, that $\unit{C^*(\G_\tb)}$ is generated by $\{1_{\X_s}: s\in
\S\} \cup \{1\}$, and contrasting \lcite{\psiPiu} with
\lcite{\PiuPsi}, we deduce that $\Lambda$ and $\widetilde\psi$ are each others
inverse.  In particular, $\Lambda$ is a *-isomorphism.

The last sentence of the statement now follows immediately from the
claim at the begining of the proof.
  \proofend

\section{Totally disconnected spaces}
  Throughout this section we will fix a Hausdorff topological space
$X$ and a basis $\Bas$ for the topology of $X$ consisting of compact
open sets.  Evidently,  in order for such a basis to exist, $X$ must 
be locally compact and totally disconnected.

Replacing $\Bas$ by the collection of all finite intersections of
members of $\Bas$, we may assume that $\Bas$ is closed under
intersections.  We will also suppose that the empty set is a member of
$\Bas$ (which is necessarily the case if $X$ has more than one point).
Summarizing, $\Bas$ is a basis for the topology of $X$ such that 
  \Bitem if $U\in\Bas$, then $U$ is compact and open,
  \Bitem $\emptyset \in \Bas$, 
  \Bitem if $U, V\in \Bas$, then $U\cap V\in \Bas$.
  \medskip

With the usual order of inclusion it is clear that $\Bas$ is a
semilattice with zero, the role of the latter being played by the
empty set.  
  This semilattice structure of $\Bas$ will be our main concern from now on
and it is our purpose to show that its tight spectrum is homeomorphic to $X$.

If $\F$ is any collection of subsets of $X$ we will adopt the notation
$\inters\F$ to refer to the intersection of all members of $\F$,
namely
  $$
  \inters\F = \bigcap_{A\in\F} A.
  \eqmark DefineIntersec
  $$
  This is in fact the standard convention for the
\stress{intersection operator} in set theory.

\state Lemma \label TopoLemma
  Let $X$ be a topological space, let $\F$ be a family of closed
subsets of $X$, and let $U$ be an open subset of $X$ such that
  $$
  \inters \F\c U.
  $$
Suppose moreover that at least one member of $\F$ is compact.  Then
there are finitely many $F_1, F_2, \ldots, F_n\in\F$, such that
$F_1 \cap F_2 \cap \ldots \cap F_n\c U$.
 
\proof Let $G\in\F$ be compact and consider the following collection
of compact subsets of $X$:
  $$
  \F' = \{F\cap G \cap (X\setminus U): F\in\F\}.
  $$
  It is clear that
  $$
  \inters\F' \c \inters\F \c U,
  $$
  but since $\inters\F'$ is also contained in $X\setminus U$, we see
that $\inters\F'=\emptyset$.  By compactness there are finitely many
$F_1, F_2, \ldots, F_n\in\F$, such that $F_1 \cap F_2 \cap \ldots \cap
F_n\cap G\cap (X\setminus U) = \emptyset$, and hence
  $$
  F_1 \cap F_2 \cap \ldots \cap F_n\cap G\c X\setminus U.
  \proofend
  $$

Given a point $x$ in $X$, the set
  $$
  \xi_x = \{U\in \Bas: U\ni x \}
  \eqmark DefineXix
  $$
  is evidently an ultra-filter, and hence also a tight filter
\ci{12.7}.  The converse of this statement is in order.

\state Proposition
  \label FoundPoint
  If $\xi$ is a tight filter on $\Bas$ then there exists a unique $x$ in
$X$ such that $\xi=\xi_x$.  In other words, the correspondence
  $$
  \Phi:  x \in X \mapsto \xi_x \in \Bast
  $$
  is a bijection.

\proof Because $\xi$ is a filter it has the finite intersection
property and, since the members of $\xi$ are compact sets, we see that
  $
  \inters \xi
  $
  (notation as in \lcite{\DefineIntersec}) is nonempty.  

Our first claim is that $\inters \xi$ consists of a single point. To
prove this assume otherwise and for each $x\in \inters\xi$, choose an
open set $W_x$ such that
  $$
  x\in W_x\not\supseteq \inters\xi,
  \subeqmark WnotContain
  $$
  (to obtain such a set one may use the Hausdorff property to separate
$x$ from any other point of $\inters\xi$).  Since $\Bas$ is a basis for
the topology of $X$ we might as well take $W_x\in\Bas$.  The $W_x$
therefore form a cover for $\inters\xi$ and, since the latter is
compact, we may choose a finite subcover, say 
  $$
  \inters\xi\c W_{x_1}\cap \ldots\cap W_{x_k}.
  $$
  Employing Lemma \lcite{\TopoLemma} we deduce that there are finitely many $F_1, \ldots, F_n\in\xi$, such that
  $$
  V:= F_1 \cap \ldots \cap F_n\c W_{x_1}\cap\ldots\cap W_{x_k}.
  $$
  Since $\xi$ is a filter,  and hence closed under finite intersections, we see that 
  $
  V
  \in\xi.
  $
  It is therefore evident that
  $$
  \{W_{x_1}\cap V, \ldots, W_{x_k}\cap V\}
  $$
  is a cover for $V$ in the technical sense of \ci{11.5}.  If $\phi$
is the character of $\Bas$ associated to $\xi$, as in \ci{12.6}, we therefore
have by \lcite{\TightCondition} that
  $$
  \bigvee_{i=1}^k\phi(W_{x_1}\cap V) = \phi(V).
  $$
  Given that $V\in\xi$, the right hand side above equals 1, and hence
there exists some $i$ such that $W_{x_i}\cap V\in\xi$, so that also
$W_{x_i}\in\xi$.  It follows that 
  $
  W_{x_i}\supseteq \inters \xi,
  $
  contradicting \lcite{\WnotContain}, and hence proving our claim.  

We conclude that $\inters\xi = \{x_0\}$, for some $x_0$ in $X$, and
the proof will be finished once we prove that $\xi=\xi_{x_0}$.  Since
the inclusion $\xi\c\xi_{x_0}$ is trivial, we concentrate on the
opposite one.  Given $U\in\xi_{x_0}$, notice that
  $$
  \inters\xi = \{x_0\} \c U,
  $$
  so we may invoke Lemma \lcite{\TopoLemma} once more to produce $F_1,
\ldots, F_n\in\xi$, such that
  $$
  W:= F_1 \cap \ldots \cap F_n\c U.
  $$
  Again we have that $W\in\xi$, so that also $U\in\xi$, concluding the
proof that $\xi_{x_0}\c\xi$.
  
The fact that $x_0$ is uniquely determined follows from the Hausdorff
property of $X$.
  \proofend

\state Proposition \label PhiUDU
  If $U\in \Bas$ then the image of\/ $U$ under the map $\Phi$ of
\lcite{\FoundPoint} coincides with the set 
  $$
  \D_ U :=
  \{\xi \in \Bast: \xi\ni U\}.
  $$

\proof For $x\in X$ we have that
  $$
  x\in U \iff
  \xi_x \ni U \iff
  \xi_x\in \D_ U \iff
  \Phi(x) \in \D_ U.
  \proofend
  $$

We may now prove the main result of this section:

\state Theorem \label TighAndX Let $X$ be a Hausdorff topological
space and let $\Bas$ be a basis for the topology of $X$ consisting of
compact open subsets,  and which is moreover closed under
finite intersections.  Assuming that $\emptyset\in\Bas$, and viewing $\Bas$
as a semilattice, the correspondence
  $$
  \Phi:  x \in X \mapsto \xi_x \in \Bast
  $$
  gives a homeomorphism from $X$ to the tight spectrum of $\Bas$.

\proof
  After  \lcite{\FoundPoint} all that remains to be proved is that
$\Phi$ is continuous and open.
A basic open set in
  $\Bast$
  is of the form 
  $$
  \Omega = \{\xi\in\Bast: U_1,\ldots,U_n\in\xi;\ V_1,\ldots,V_m\not\in\xi\},
  $$
  where $U_1,\ldots,U_n;\ V_1,\ldots,V_m\in\Bas$.  With notation as in
\lcite{\PhiUDU} notice that
  $$
  \Omega =   \bigcap_{i=1}^n \D_{U_i}\ \cap\ \bigcap_{j=1}^m \Bast\setminus \D_{V_i},
  $$
  so by \lcite{\PhiUDU} we have
  $$
  \Phi\inv(\Omega) = \bigcap_{i=1}^n U_i\ \cap\ \bigcap_{j=1}^m
X\setminus V_i,
  $$
  which is an open subset of $X$, proving the continuity of $\Phi$.
To prove that $\Phi$ is an open map it suffices to verify that
$\Phi(U)$ is open in $\Bast$, for all $U\in \Bas$, but this follows
at once from \lcite{\PhiUDU}.  \proofend
  
\section{Groupoids with totally disconnected unit space}
  In this section we will fix a (not necessarily Hausdorff) \'etale
groupoid $\G$, with range map denoted $\r$, and domain (or source) map
written $\d$.  For a precise definition see e.g., \cite{\Paterson} or
\scite{\Actions}{3.1}.  

Recall that a \stress{{\slice}} of $\G$ is any open subset $S\c \G$,
such that $\r$ and $\d$ are injective when restricted to $S$.  It is
well known,  and easy to prove,  that the collection of all {\slice}s forms
a basis for the topology of $\G$.

If $S$ and $T$ are {\slice}s, their product
  $$
  ST = \{\sigma\tau: \sigma\in S,\ \tau\in T,\ \d(\sigma) = \r(\tau)\},
  $$
  may be shown to be a {\slice}.  With this multiplication operation
the collection 
of all {\slice}s  forms a semigroup, and in fact an inverse
semigroup, where the
\stress{inverse} of a {\slice} $S$ is given by
  $$
  S^* = \{\sigma\inv: \sigma\in S\}.
  $$

  By trivial reasons every open subset of $\Gz$ (the unit space of $\G$) is a
{\slice}, so the inverse semigroup of all {\slice}s contains the
topology of $\Gz$.  The latter is in fact  identical to the
idempotent semilattice of the former.

Recall that a topological space is said to be \stress{totally
disconnected} when its topology admits a basis formed by 
subsets which are both open and closed.

\state Proposition \label basisOfCompactBissections If $\Gz$ is totally
disconnected then the collection of all compact {\slice}s forms a
basis for the topology of $\G$.

\proof
  Given $\g\in U\c \G$, with $U$ open, we need to find a compact
{\slice} $T$ such that 
  $$
  \g\in T \c U.
  \subeqmark ThisTMustSatisfy
  $$
  Choose a {\slice} $S$ such that $\g\in S\c U$, and notice that
since
  $
  \d(\g) \in \d(S),
  $
  and since $\d(S)$ is open, there exists a compact\fn{Notice that
the axioms for \'etale groupoids include the requirements that $\Gz$ be
Hausdorff and 
locally compact.  With the added assumption that $\Gz$ is totally
disconnected one may easily prove that the collection of all open
\stress{compact} subsets is a basis for the topology of $\Gz$.} open subset
$K$ of $\Gz$ such that
  $
  \d(\g)\in K\c \d(S).
  $
  Observing that $\d$ is a homeomorphism from $S$ to
$\d(S)$, we see that $T := S\cap d\inv(K)$ is a compact
{\slice} satisfying \lcite{\ThisTMustSatisfy}.  \proofend

Since we are not assuming $\G$ to be Hausdorff, the compact slices
provided by the above result might not be closed subsets.  We have
therefore fallen short of proving that $\G$ is totally disconnected,
at least according to the classical definition of this concept which
requires a basis consisting of subsets which are both open and
\stress{closed}.  I wonder if the definition of totally disconnected
non Hausdorff spaces should be adapted to include spaces such as our
$\G$ above.

It is easy to prove that the product of two compact {\slice}s is
also compact,  so the collection of all compact {\slice}s
forms a semigroup.  It is clearly also an inverse semigroup since
$S^*$ is easily seen to be compact, should $S$ be compact.

In order to proceed we need the following elementary result in General
Topology, which is perhaps well known to the reader:

\state Lemma A (not necessarily Hausdorff) second countable 
topological space contains at most countably many compact open
subsets.

\proof Let $\Bas$ be a countable basis for the topology
of our space.  Given a compact open subset  $K$, for each $x\in K$,
pick $U_x\in \Bas$ such that $x\in U_x\c K$.  Thus $\{U_x\}_{x\in K}$ is
an open cover for $K$, from which we may extract a finite subcover,
say
  $\{U_{x_i}\}_{i=1}^n$.
  Consequently $K$ is the union of finitely many members of $\Bas$.  Since
there are only countably many finite subsets of $\Bas$, the result follows.
  \proofend

If we assume that, besides satisfying the hypothesis of
\lcite{\basisOfCompactBissections}, our groupoid is second countable,
it will follow that the inverse semigroup of all compact
{\slice}s is countable.

Summarizing, we have:

\state Corollary Let $\G$ be a (not necessarily Hausdorff) \'etale,
second countable groupoid with totally disconnected unit space.
Then there is a (necessarily) countable inverse semigroup $\S$ consisting of
compact {\slice}s of $\G$ which is moreover a basis for the topology  of $\G$.

\proof
  The collection $\S$ formed by \stress{all} compact {\slice}s
fits the above description. \proofend

From now on we fix a countable inverse semigroup $\S$ consisting of
compact {\slice}s of $\G$ which is a basis for the topology  of
$\G$, and let us denote by $\E$ \ the idempotent semilattice
of $\S$.  
Each member of $\E$ is therefore a compact open subset of $\Gz$.  
It is also clear that $\E$ is a basis for the topology of $\Gz$,  so   by
\lcite{\TighAndX} we obtain a homeomorphism 
  $$
  \Phi: \Gz\to \Et,
  \eqmark PhiForGzero
  $$
  such that $\Phi(U) = \D_ U$, for all $U\in \E$, according to \lcite{\PhiUDU}.


Recall from \ci{5.3} that $\S$ acts on $\Gz$ as follows:
  given $S\in\S$,
  $$
  \lambda_S: \d(S)\to\r(S)
  $$
  is defined to be the homeomorphism given,
for every $x\in \d(S)$, by
  $$
  \lambda_S(x) = \r(\gamma),
  $$
  where $\gamma$ is the unique element in $S$ such that
$\d(\gamma)=x$.  This action was denoted by $\t$ in \cite{\Actions}
but we denote it here by $\lambda$ in order to distinguish it from
the action of $\S$ on its tight spectrum,  as in \lcite{\DefineTheta}.

\state Lemma \label DgSUSrgU Let $S$ and $U$ be {\slice}s of $\G$
with $U\c\Gz$. Then, for every $\g\in S$, one has that
  $\d(\gamma)\in S^*US$, if and only if  $\r(\gamma) \in U$.

\proof Assuming that $\d(\gamma)\in S^*US$, write $\d(\gamma) =
\alpha\inv y \beta$, with $\alpha,\beta\in S$,  and $y\in U$.
Then 
  $
  \d(\g) = \d(\b),
  $
  and hence $\g=\b$, because $\d$ is injective on $S$.  It follows that
  $$
  \r(\g) = \r(\b) = \d(y) = y \in U.
  $$
  Conversely, if $\r(\g) \in U$, then
  $$
  \d(\g) = \g\inv \g = \g\inv \r(\g) \g \in S^*US.
  \proofend
  $$

\state Proposition The homeomorphism $\Phi$ in \lcite{\PhiForGzero}
establishes an equivalence between $\lambda$ and $\actspec$.  More
precisely, for every $S\in\S$, the diagram below commutes:
  $$
  \matrix{
  \d(S) & \buildrel \lambda_S \over \longrightarrow & \r(S)\cr\cr
  \Phi \downarrow \ \ &&\ \ \downarrow \Phi\cr\cr
  \D_{\d(S)} & \buildrel \actspec_S \over \longrightarrow & \D_{\r(S)}
  }
  $$

\proof  Let $x\in\d(S)$, and write $x = \d(\g)$, with $\g\in S$.
Then, using brackets to denote boolean value,  for every $U\in \E$, one has
  $$
  \actspec_S(\Phi(x))\calcat U = 
  \Phi(x)\calcat{S^*US} =
  [\d(\g)\in S^*US] \={(\DgSUSrgU)}
  [\r(\g)\in U] \$=
  \Phi(\r(\g))\calcat U =
  \Phi(\lambda_S(x))\calcat U.
  \proofend
  $$

It is evident that the equivalent actions of $\S$ above give rise to
isomorphic groupoid of germs, so the groupoid of germs for the action
$\lambda$, which we will temporarily denote by $\H$, is isomorphic to the groupoid
of germs for the action $\actspec$, which we have denoted by
$\G_\tightbox$ in section \lcite{\PrelimSection}.  The obvious isomorphism, say
  $f: \H \to \G_\tightbox$,
  is then given by
  $f([S,x]) = [S,\Phi(x)]$,
  for all $S\in\S$, and $x\in \d(S)$.

In addition, recall from \ci{5.4} that $\H$ is isomorphic to $\G$,
under an isomorphism 
  $g:\H\to \G$
  which sends the germ $[S,x]$
  to the unique element $\gamma\in S$ such that $\d(\gamma) = x$.  We
thus get an isomorphism from $\G_\tightbox$ to $\G$ by composition,
namely
  $$
  \G_\tightbox \buildrel f\inv \over \to \H \buildrel g \over \to \G.
  \eqmark IsoGpds
  $$

This shows that $\G$ may be recovered from $\S$, a result which
deserves to be highlighted in the following:

\state Theorem Let $\G$ be a (not necessarily Hausdorff) \'etale,
second countable groupoid with totally disconnected unit space.  Also
let $\S$ be a (necessarily) countable inverse semigroup consisting of
compact {\slice}s of $\G$ which is a basis for the topology of $\G$.
Then $\G$ is isomorphic to the groupoid 
$\G_\tb$ formed by the germs of the canonical action of $\S$ on the
tight spectrum of its idempotent semilattice.

Given $S\in\S$, let $\X_S\c\G_\tb$ be given by
  $$
  \X_S = \{[S,\xi]: \xi\in \D_{\d(S)}\},
  $$
  and recall that the characteristic function $1_{\X_S}$ was referred
to in \lcite{\Reinterpret} as being equal to $\Lambda(\piu(S))$.  We
are interested in determining the image of $\X_S$ under $g\circ f\inv$. 
  Clearly
  $
  f\inv(\X_S) = \{[S,x]: x\in d(S)\},
  $
  whence
  $
  g\big(f\inv(\X_S)\big) = S.
  $
  Thus,   if we let 
  $$
  \Gamma:   C^*(\G_\tightbox) \to C^*(\G).
  $$
  be the *-isomorphism induced from \lcite{\IsoGpds}, we conclude from \lcite{\Reinterpret} that 
  $$
  \Gamma(\Lambda(\piu(S))) = \Gamma(1_{\X_s}) = 1_S.
  $$

We have therefore proven the following, which is the main result of
this paper:

\state Theorem Let $\G$ be a (not necessarily Hausdorff) \'etale,
second countable groupoid with totally disconnected unit space.  Also
let $\S$ be a (necessarily) countable inverse semigroup consisting of
compact {\slice}s of $\G$ which is a basis for the topology of $\G$.
Then there is a *-isomorphism
  $$
  \psi :   C^*_\tightbox(S) \to  \unit{C^*(\G)},
  $$
  such that $\psi\circ\piu(S) = 1_S$,  for all $S\in \S$.

\proof It is enough to take $\psi = \widetilde\Gamma\circ\Lambda$,
where $\widetilde\Gamma$ is the extension of $\Gamma$ to the
corresponding unitized algebras.  \proofend

Borrowing from the universal property of $C^*_\tightbox(S)$ we have the
following main consequence:

\state Corollary
  Let $\G$ and $\S$ be as above.  Then for every tight representation
$\pi$ of $\S$ in a unital C*-algebra $A$, there exists a unique unital
*-homomorphism $\phi:\unit{C^*(\G)} \to A$, such that
  $\phi(1_S) = \pi(S)$,  for every $S\in\S$.


  \endgroup

  \begingroup
  \bigskip\bigskip 
  \font \sc = cmcsc8 \sc
  \parskip = -1pt

  Departamento de Matem\'atica 

  Universidade Federal de Santa Catarina

  88040-900 -- Florian\'opolis -- Brasil

  \eightrm r@exel.com.br

  \endgroup
  \end